\documentclass[12pt]{amsart}
\usepackage{amsthm}
\usepackage{amssymb}
\usepackage{amsmath}
\usepackage{amscd}
\usepackage{epic}
\usepackage{eepic}

\newcommand{\printname}[1]
  {\smash{\makebox[0pt]{\hspace{-2.0in}\raisebox{8pt}{\tiny #1}}}}
\newcommand{\labell}[1] {\label{#1}\printname{#1}}

\newtheorem {Theorem}   {Theorem}
\newtheorem {Lemma} {Lemma} [section]

\theoremstyle{definition}

\theoremstyle{remark}
\newtheorem{Remark}[Lemma]{Remark}

\newtheorem {Corollary}[Lemma]{Corollary}
\newcommand {\cal}[1]  {{\mathcal{#1}}}

\newcommand{\nab}[1][]{\ensuremath{\mathrm{\nabla}{#1}}}

\newcommand{\be}{\begin{equation}}
\newcommand{\ee}{\end{equation}}

\newcommand{\del}{\partial}

\newcommand{\lb}[1]{\label{#1}}
\newcommand{\lbl}[1]{\labell{#1}}

\newcommand{\De}{\Delta}

\newcommand{\al}{\alpha}
\newcommand{\om}{\omega}
\newcommand{\Om}{\Omega}
\newcommand{\we}{\wedge}

\newcommand{\sig}{\sigma}

\newcommand{\Ref}[1]{(\ref{#1})}
\newcommand{\fr}{\frac}

\newcommand{\ra}{\rightarrow}

\def\t{\tau}

\newcommand{\wht}{\widehat}

\begin{document}

\title[The eta invariant in the doubly K\"ahlerian conformally compact\ldots]
{The eta invariant in the doubly K\"ahlerian conformally compact Einstein case}

\author[]{}
\address{}
\date{\today}

\author[Gideon Maschler]{Gideon Maschler}
\address{Department of Mathematics and Computer Science, Clark University,
Worcester, Massachusetts 01610, U.S.A.}
\email{gmaschler@clarku.edu}

\thanks{}

\begin{abstract}
On a $3$-manifold bounding a compact $4$-manifold, let a conformal structure be induced from a
complete Einstein metric which conformally compactifies to a K\"ahler metric. Formulas
are derived for the eta invariant of this conformal structure under additional
assumptions. One such assumption is that the K\"ahler metric admits a special
K\"ahler-Ricci potential in the sense defined by Derdzinski and Maschler. Another
is that the K\"ahler metric is part of an ambitoric structure, in the sense defined
by Apostolov, Calderbank and Gauduchon, as well as a toric one. The formulas are derived using the
Duistermaat-Heckman theorem. This result is closely related to earlier work of Hitchin
on the Einstein selfdual case.
\end{abstract}

\maketitle

\section{Introduction}
In \cite{h}, Hitchin computed an obstruction, in terms of a bound on the eta invariant for
a conformal structure on the $3$-sphere to be induced from a complete self-dual Einstein metric
on the $4$-ball. Following LeBrun \cite{l}, such conformal structures are said to have
\emph{positive frequency}, drawing on an analogy between this case and the classical obstruction
for a smooth function on the circle to be the boundary value of a holomorphic function on the disk.

In this paper we give formulas for the eta invariant of a conformal structure induced from another
type of asymptotically hyperbolic Einstein metric in dimension four, namely one that is conformal to a
(conformally compact) K\"ahler metric, and satisfies additional technical assumptions. These formulas
are obtained via moment map techniques, specifically the Duistermaat-Heckman Theorem.
The assumptions are either that the K\"ahler metric admits a special K\"ahler-Ricci potential in the sense
of \cite{dr-ma1,dr-ma2}, or that it is part of an ambitoric structure, in the sense defined recently
by Apostolov et. al. \cite{acg}, as well as a toric one. Common to both cases is the existence of a second,
oppositely oriented complex structure, with respect to which one has a second K\"ahler metric in the
conformal class, which plays a major role in the derivation of the formulas.

Section \ref{two} is devoted mainly to a review of the theorem of Hitchin. Section \ref{skrbd} gives
an overview of the classification of metrics with a special K\"ahler-Ricci potential in the case of a manifold
with boundary, and also discusses the dual K\"ahler metric. In Section \ref{formulas} the two formulas
are derived.




\section{Preliminaries}\lb{two}
\subsection{Conformal compactifications of Einstein metrics}
In the following let $(M,g_E)$ be \emph{conformally compact Einstein}, i.e.,
$M$ is a compact manifold with a nonempty boundary, $g_E$ a complete Einstein metric in the interior
of $M$ and there exists a smooth defining function $\tau$ for the boundary
$\del M$  (so that $\tau\geq 0, \tau|_{\del M}=0, d\tau|_{\del M}\neq 0$) for which
$$\text{$g=\tau^2g_E$ is smooth on $M$.}$$ The pair $(g,\tau)$ is called a {\em
conformal compactification} for $g_E$. The restriction $\gamma=g|_{\del M}$ varies
with $\tau$ within a fixed conformal class $[\gamma]$ on $\del M$.

Later we will consider such $(M,g_E)$ with a \emph{K\"ahlerian (conformal) compactification}, i.e.
satisfying the additional assumption that $g$ is K\"ahler in the interior of $M$.
The K\"ahlerian assumption, combined with the more technical assumption \Ref{skr}, which will be
described in \S\ref{skrbd}, imply that the function $Q:=|\nab\t|^2$
is constant on each level set of $\t$ (see \cite{dr-ma1}). A compactification with the latter property
may be called {\em uniform}. It is easier to achieve, in a sense, than the more well-known
{\em geodesic} one, where $Q$ is constant in a collar neighborhood of the boundary: from a given
compactification, one must solve a PDE to obtain a geodesic compactification, but only an ODE to
obtain a uniform one. On the other hand, the above two assumptions imply the uniformity of the compactification
throughout the open and dense set of non-critical $\tau$ values, and not just in a collar neighborhood of
the boundary.

\subsection{The theorem of Hitchin}
We recall here the result of Hitchin, emphasizing its applicability, in addition to the four-ball,
also to other four-manifolds.
\begin{Theorem}[Hitchin]\lb{hit}
Suppose $(g,\tau)$ is a conformal compactification of $(M,g_E)$, 
with $g_E$ Einstein and $M$ a manifold with boundary of real dimension four. Then
\be\lb{hitchin}\eta_\gamma(\del M)=\fr1{12\pi^2}\int_M\left(|W_-|^2_g-|W_+|^2_g\right)
\text{vol}_{g}-\sig(M),\end{equation}
where $W_+$ ($W_-$) are the (anti-)self dual parts of the Weyl tensor of $g$, $\text{vol}_g$
its volume form, $\sig(M)$ the signature of $M$ and $\eta_\gamma(\del M)$ the eta invariant of the
boundary metric $\gamma$.
\end{Theorem}
Here the eta invariant is defined as the value at $s=0$ of the series
$\sum_{\lambda\neq 0}(\text{sgn}\lambda)|\lambda|^{-s}$ (which is holomorphic for $\text{Re}\,s>-1/2$),
where the summation
is over nonzero eigenvalues of the self-adjoint operator on even forms of $\del M$, given by
$$B(\al)=(-1)^p(*d-d*)\al, \quad \al\in\Omega^{2p}(\del M).$$ This invariant depends only on the conformal
class of $\gamma$ in $\del M$.
\begin{proof}
The Atiyah, Patodi and Singer \cite{aps} signature formula for a manifold with boundary reads \cite{egh}
$$\sig(M)=-\fr 1{24\pi^2}\int_M\,\mathrm{tr}(R\we R)+\fr 1{24\pi^2}\int_{\del M}\text{tr}(\Pi\we R)
-\eta_\gamma(\del M),$$
where $R$ is the curvature tensor of $M$ and $\Pi$ is the second fundamental form of $\del M$, considered
as a one-form valued in (restrictions to $\del M$ of) endomorphisms of the tangent bundle of $M$.
In the case at hand, this boundary term actually vanishes. In fact, as $g$ is conformal to the Einstein metric $g_E$,
it satisfies a Ricci-Hessian equation of the form $\nab d\tau + (\tau/2)\mathrm{r}=\sigma g$, with $\mathrm{r}$ the
Ricci curvature and $\nabla$ the covariant derivative of $g$, for some function $\sigma$. As $\del M=\{\tau=0\}$,
the second fundamental form of the boundary is the Hessian $\nab d\tau$. But the Ricci-Hessian equation above 
yields that on the boundary this Hessian is a multiple of the metric, i.e. that $\del M$ is totally umbilical. 
Thus, in an adapted orthonormal frame $\{e_i\}_{i=1}^4$ with $\{e_i\}_{i=1}^3$ tangent to, and $e_4$ a unit normal 
to $\del M$, $\Pi=k\sum_{i=1}^3e_i\otimes e_i$ for some constant $k$, which implies
$\text{tr}(\Pi\we R)=k(e_1\we R_{14}+e_2\we R_{24}+e_3\we R_{34})+k(R_{1423}+R_{2431}+R_{3412})e_1\we e_2\we e_3)
=k(R_{4123}+R_{4231}+R_{4312})e_1\we e_2\we e_3)=0$ by the Bianchi identity. The result now follows as
$\mathrm{tr}(R\we R)=2\left(|W_+|^2_g-|W_-|^2_g\right)\text{vol}_{g}$.
\end{proof}

For the four-ball, $\sig(M)=0$, and so Hitchin concludes from \Ref{hitchin} that the eta invariant of $S^3$
is nonpositive for conformal structures induced from self-dual Einstein metrics on $B^4$, and zero exactly
for the standard structure. More generally, the following obvious general bounds hold:
\begin{Corollary}
\be\lb{gen}-\fr1{12\pi^2}\int_M|W_+|^2_g\,\text{vol}_{g}-\sig(M)\leq\eta_\gamma (\del M)\leq \fr1{12\pi^2}\int_M|W_-|^2_g\,\text{vol}_{g}-\sig(M)\end{equation}
with equality holding on the left-hand side in the self-dual case, and on the right-hand side in
the anti-selfdual case. Moreover, if the compactification is also Kahlerian, and $g$ is not anti-selfdual,
then the left hand inequality of \Ref{gen} can be rewritten
as \be\lb{scalar}-\fr1{288\pi^2}\int_M s_g^2\text{vol}_g-\sig(M)
\leq \eta_\gamma (\del M),\end{equation}
an inequality that holds because on a K\"ahler surface, $W_+$ is completely determined 
by the scalar curvature and the K\"ahler form (see for example Derdzinski \cite{der}).
\end{Corollary}

\section{K\"ahler-Ricci potentials on manifolds with boundary}\lb{skrbd}
\subsection{Metrics with a K\"ahler-Ricci potential}\lb{outline}

We give here an outline of the analog, for manifolds with boundary, of the main theorem in \cite{dr-ma2}
(see Theorem \ref{bihol-iso} below).

Let $(M,J,g,\t)$ be a quadruple consisting of a compact manifold $M$ with a nonempty boundary,
 having a defining function $\t$, and a complex structure
$J$ along with a K\"ahler metric $g$ on the interior of $M$ for which
\be\lb{skr} \text{$\t$ is a special K\"ahler-Ricci potential.}\end{equation}
This requirement on $\t$ is a technical notion which holds in the conformally Einstein case if either the
complex dimension of $M$ is at least three, or else it is two and
$d\t\wedge d\De_g\t=0$, with $\De_g$ the $g$-Laplacian. Metrics satisfying \Ref{skr} have been classified
both locally and on compact manifolds in \cite{dr-ma1} and \cite{dr-ma2}. Specifically, a special K\"ahler-Ricci
potential is a Killing potential on a K\"ahler manifold, i.e. $J\nab\t$ is a Killing vector field,
for which the following additional condition holds: at any noncritical point of $\t$,
the subspace ${\cal{H}}=(\mathrm{span}(\nab\t, J\nab\t))^\perp$ of the tangent space is an eigenspace
for both the Hessian $\nab d\t$ and the Ricci curvature $\mathrm{r}$.

Suppose the interval of $\t$-values is $[0,\t_0]$ and the $\t$-preimage $N$ of the maximal value $\t_0$ is a critical manifold.
By \cite[Lemma 7.3]{dr-ma2}, $N$ is a complex submanifold of $M$ which is of complex codimension one, or else consists of a single point.
Furthermore, {\em any} critical manifold is the $\t$-preimage of an extreme value for $\t$ (see the proof of Proposition $11.4$
in \cite{dr-ma2}). As $\t$ is a defining function, the zero level set of $\t$ does not contain any critical points (and is the boundary),
so that $N$ is the unique critical manifold of $\t$.

Special K\"ahler-Ricci potentials give rise to ingredients that are used to give an explicit construction of a metric $\overline{g}$
which is biholomorphically isometric to $g$.
Namely, the function $|\nab\t|^2$ is a composite consisting of $\t$ followed by a smooth function $Q(\t)$ defined on $[0,\t_0]$.
This follows from Lemma 9.1 of \cite{dr-ma2} along with Proposition 11.5(ii)  there, for the interval $[0,\t_0]$
with the proof changed to the case of a compact manifold with boundary. Adopting again the
same proposition, $Q$ satisfies the boundary and positivity conditions
\be\lbl{bdps} \text{$Q(\t_0)=0$, $Q'(\t_0)\ne 0$ and $Q(\t)>0$ for $\t\in [0,\t_0)$.}\end{equation}
Another ingredient is a constant $c$, defined as follows: the eigenfunction $\phi$ of the Hessian $\nab d\t$
on the subspace ${\cal{H}}$ can be identically zero, in which case $c$ is undefined, or nowhere zero
for $\t$-noncritical points, and then $c=\t-Q/(2\phi)$ \cite[Lemma 3.1]{dr-ma2}. We note for future reference that \be\lbl{ct0}
\text{$c\notin [0,\t_0]$ unless
the critical manifold $N$ consists of a single point and then $c=\t_0$}\end{equation}
(see \cite[Lemma 7.5]{dr-ma2}). Furthermore, the unique nonzero eigenvalue of $\nab d\t$
at each point of $N$ is a constant denoted $e$ \cite[Proposition 7.3 (denoted $a$ there)]{dr-ma2}. With $Q,c$ and $e$ one
defines a metric $\overline{g}$ on a disc bundle $S$ contained in the total space of the normal bundle ${\cal{L}}$
to $N$, using a Hermitian fiber metric $\langle\cdot,\cdot\rangle$ whose real part is $g|_{\cal{L}}$ and a metric $h$ on $N$
defined as either $g|_{TN}$ if $\phi=0$, or $(2|\t_0-c|)^{-1}g|_{TN}$ if $\phi \ne 0$ (if $N$ consists of a single point
$\{y\}$ then $S$ is a disc, ${\cal{L}}=T_yM$ and $h$ is undefined). The connection associated with the fiber metric is the normal
connection, and its curvature $\Omega$ satisfies \be\lb{curv-L} \Om=p\,\omega^h,\end{equation} where $\omega^h$ is
the K\"ahler form of $h$ and the constant $p$ equals $\pm 2e$ (see \cite[Lemma 12.4]{dr-ma2}).

Regarding $N$ as the zero section $N\subset {\cal{L}}$, the metric $\overline{g}$ is defined as follows:
\be\lbl{met}\end{equation}
\begin{tabbing}
\text{\emph {if\ $N$\ is\ a\ complex\ codimension}}\=\text{\emph{\ one\ submanifold:\ }}\=\nonumber\\[5pt]
\>i] $\overline{g}|_{\cal{H}}=\pi^*h$ \>\quad\text{\ if\ }$\phi=0$,\nonumber\\[4pt]
\>ii] $\overline{g}|_{\cal{H}}=2|\t-c|\,\pi^*h$ \>\quad\text{\ if\ }$\phi\ne 0$,\nonumber \\[4pt]
\text{\emph {if\ }}$N$\ \text{\emph{consists\ of\ a\ single\ point:\ }}\>\nonumber\\[4pt]
\>iii] $\overline{g}|_{\cal{H}}=(2|\t-c|/(|e|r^2))\,\pi^*h,$\nonumber\\[4pt]
\text{\emph {and\ in\ all\ cases:\ }}\>\nonumber\\[4pt]
\>$\overline{g}|_{\cal{V}}=\left( Q(\t)/(er)^2\right)\,
\mathrm{Re}\,\langle\cdot,\cdot\rangle.$
\end{tabbing}
Here, if $N$ is a complex codimension one submanifold, ${\cal{H}}, {\cal{V}}$ are the horizontal distribution defined by the
connection and the vertical distribution, respectively, which are also declared $\overline{g}$-orthogonal to each other. If $N$ consists
of a single point, then ${\cal{V}}$ is the distribution on the vector space $T_y{L}\setminus 0$ defined at each nonzero vector as the
complex line through this vector, while ${\cal{H}}$ is the distribution orthogonal to ${\cal{V}}$ with respect to
$\mathrm{Re}\,\langle\cdot,\cdot\rangle$. Finally, $r=r(\t)$ is a nonnegative function on $[0,\t_0]$ which vanishes only at
$\t=\t_0$ and satisfies the differential equation $dr/d\t=e\,r/Q$ in $(0,\t_0)$. Abusing notation by regarding $r$ also
as the norm function of $\langle\cdot,\cdot\rangle$, it is a function on ${\cal{L}}$. Using the inverse map
$\t(r)$ to $r(\t)$, all functions of $\t$ in \Ref{met} become functions of the norm $r$ on ${\cal{L}}$,
thus giving rise to a metric defined on $\cal{L}$. Furthermore, there exists a positive value $r_0$
for which $\t(r_0)=0$, and the disc bundle $S$ is characterized by the inequality $r\leq r_0$.

A version of the main theorem of \cite{dr-ma2}, Theorem 16.3 there, now holds for $M$ with essentially the same proof, i.e.
\begin{Theorem}\lb{bihol-iso}
Let $M$ be a compact manifold with a nonempty boundary, with a K\"ahler metric $g$
and a special K\"ahler-Ricci potential $\t$ on its interior, both extending smoothly to the boundary. 
Suppose $\t$ has the value zero on the boundary. Then the triple $(M,g,\t)$ is biholomorphically
isometric to the disc bundle (or disc) $(S,\overline{g},\t(r))$ described above, and the orbits of
the Killing vector field $J\nab \t$ are closed.
\end{Theorem}
Note that the assumption on the boundary value of $\t$ excludes the possibility of $S$ being
an annulus bundle, and guarantees that $M$ has a critical manifold, the preimage of the nonzero extremum
$\t_0$ of $\t$. Also, the normal exponential map which is a building block for the inverse of
the biholomorphic isometry gives a diffeomorphism $S\setminus \{r=r_0\}\ra M\setminus \del M$
(this is \cite[Lemma 13.2]{dr-ma2}, adjusted to the case of a manifold with boundary).


The last part of the conclusion of the theorem, which can be considered a consequence of the fact
that the orbits of $J\nab\t$ coincide under the above biholomorphic isometry with the $S^1$ orbits of the
natural circle action on $L$, can be seen directly as follows. At a critical point of a special 
K\"ahler-Ricci potential $\t$, the Hessian $\nab d\t$ has exactly one nonzero eigenvalue 
(see \cite[Proposition 7.3]{dr-ma2}). By \cite[Lemma 10.2]{dr-ma2}, the vector field $J\nab\t$,
being the image under the differential of the exponential map of a linear vector field with a periodic flow,
has itself a periodic flow in a neighborhood of the critical point. By the unique continuation
property for isometries this periodicity is in fact global, as two isometries in the flow that agree 
on a nonempty open set actually agree globally. However, comparing isometries of a flow, and, more specifically,
verifying the existence of a circle action produced by a flow \cite[Corollary 10.3]{dr-ma2}, may be performed
if the vector field is complete. In \cite{dr-ma2} completeness of the metric, which
implies completeness of the vector field, was assumed. However only the latter is necessary, and
in the setting at hand, of a manifold with boundary, it occurs automatically. In fact, as $g$ and $\t$ are 
smooth everywhere, $\nab\t$ and the Killing field $J\nab\t$ are well-defined on the boundary, and since the 
latter is orthogonal to to the former, and $\nab\t$ is normal to the boundary, $J\nab\t$ is in fact tangent 
to the boundary at boundary points, and hence necessarily complete. (Note as an aside that a chain of 
reasoning exists which, starting with the fact that $J\nab\t$ is Killing, concludes with smoothness of $\t$.
Thus in the statement of Theorem \ref{bihol-iso}, it is enough to assume continuity of $\t$ at the boundary.)

\begin{Remark}
If we now add to the above assumptions the requirement on $(M,g)$ that it is a (K\"ahlerian) conformal
compactification of an {\em Einstein} metric on the interior and $M$ is of dimension four, the metric
$\overline{g}$ (and hence indirectly $g$) can be specified by a local formula for $Q$ \cite[\S19-21]{dr-ma1}. In fact, $Q(\t)$ is
a member of one of the following three families of rational functions, whose formulas are given below with
$K$, $\alpha$, $\beta$, $A$, $B$ and $C$ denoting constants (the precise choice of which ensures that $Q(\t)$
satisfies the boundary and positivity conditions \Ref{bdps}): 
\begin{eqnarray*}
\text{a]   }Q&=&-K\tau^2+(\alpha\tau^3-\beta/2)/3 \text{ if $\phi=0$.}\\
\text{b]   }Q&=&-K\tau/2+\alpha\tau^3-\beta/3 \text{ if $\phi\ne 0$ and $c=0$.}\\
\text{c]   }Q&=&({\tau/ c}-1)(AE({\tau/ c})+BF({\tau/ c})+C)\text{ if $\phi\ne 0$ and $c\ne 0$, where}\\
&&E(x)=x^2-1,\ \
F(x)=(x-2)x^3/(x-1)^2.
\end{eqnarray*}
These formulas are given for completeness, and will not be used again.
\end{Remark}

\subsection{The dual K\"ahler metric}

As in \S\ref{outline}, Let $(M,g,\t)$ be a triple consisting of a compact manifold $M$ with a nonempty boundary,
having a defining function $\t$, and a K\"ahler metric $g$ on the interior of $M$ for which $\t$ is a special
K\"ahler-Ricci potential. 
Suppose also that the isometric metric $\overline{g}$ of the last section is either
of type ii] or type iii] in \Ref{met}. Such a metric was called \emph{nontrivial SKR} in \cite{m}. In this case
the constant $c$ is defined (see paragraph between \Ref{bdps} and \Ref{ct0}) and in fact
\be\lbl{dual}\text{$\wht{g}=g/(\t-c)^2$ is K\"ahler for an oppositely
oriented complex structure}\end{equation}
while $\wht{t}=1/(\t-c)$ is a special K\"ahler-Ricci potential for $\wht{g}$ (see \cite[Proposition 5.1]{m}
or \cite[Remark 28.4]{dr-ma3}). Our main use of this fact is in the formula
\be\lbl{hat} |W_-|^2=(\t-c)^{-4}|\wht{W}_+|^2_{\wht{g}}\end{equation}
which follows since $W_-$ is, as a $(3,1)$ tensor, a conformal invariant (and using the conformal factor
$(\t-c)^2$ relating $g$ with $\wht{g}$).



\section{Moment map formulas for the $\eta$ invariant}\lb{formulas}
\subsection{The case of a special K\"ahler-Ricci potential}

In the following theorem, the meaning of the notations $N$, $h$, $c$, $e$, $\t_0$, $r_0$ and $p$ is that given in \S\ref{outline},
and Theorem \ref{bihol-iso} is implicitly used in its proof.
Thus we associate with the K\"ahler metric $g$ data used
in the construction of its isometric metric $\overline{g}$. Additionally, we continue to denote the signature
of a manifold by $\sigma(M)$, while  the volume of $N$ with respect to $h$ will be denoted $Vol_h(N)$.
\begin{Theorem}\lb{s-s}
Let $M^4$ be a compact four-manifold having a nonempty boundary with a complete Einstein metric $g_E$ on its interior,
admitting a K\"ahlerian conformal compactification $(g,\t)$ satisfying \Ref{skr} on its interior for a nontrivial SKR metric.
Then the eta invariant of the boundary metric $\gamma=g|_{\del M}$ is given by
\begin{multline}\lb{eta-short}
\eta_\gamma(\del M)=\fr 1{288\pi^2}Vol_h(N)\int_{0}^{\t_0} \left[\left(b^2(t-c)^{-6}-a^2\right)t^2\right]
\cdot \left[l+p\,t\,\right]\,dt-\sig(M)
\end{multline}
for some nonzero constants $a$ and $b$,
with $l$ standing for the constant $2|c|$ if $\t$ has
a complex codimension one critical manifold, and $2|c|/(|e|r_0^2)=2|\t_0|/(|e|r_0^2)$ (see \Ref{ct0})
if the $\t$-critical manifold consists of one point.
\end{Theorem}
\begin{Remark}
The integral in \Ref{eta-short}, is of course, elementary and equals
$$\left.\left[-\,\frac {{t}^{2}{a}^{2}p}2-
{a}^{2}lt-\,{\frac {{b}^{2}p}{ 4\left( t-c
 \right) ^{4}}}-\,{\frac {{b}^{2} \left( l+cp \right) }{ 5\left( t-c
 \right) ^{5}}}\right]\right|_0^{\t_0}.$$
\end{Remark}
\begin{proof}
By Theorem \ref{hit}, formula \Ref{hitchin} gives an expression for the eta invariant in terms of the signature, the
Weyl tensors $W_+$, $W_-$ and the volume form of $g$. Using the dual metric $\wht{g}$, one can replace $|W_-|^2$ by
expression \Ref{hat} involving $|\wht{W}_+|^2_{\wht{g}}$. Since both $g$ and $\wht{g}$ are K\"ahler metrics, the
self-dual part of their Weyl tensors can be expressed in terms of the scalar curvatures $s$ and $\wht{s}$
(as in \Ref{scalar}, cf. \cite{der}). We thus have
 \begin{multline}\lb{eta-main}
\eta_\gamma(\del M)=\fr 1{12\pi^2}\int_M \left[(\t-c)^{-4}|\wht{W}_+|_{\wht{g}}^2-|W_+|^2\right]\text{vol}_g-\sig(M)\\
=\fr1{288\pi^2}\int_M \left[(\t-c)^{-4}\,\wht{s}^{\ 2}-s^2\right]\text{vol}_g-\sig(M)
\end{multline}
Next, in dimension four, the scalar curvature of a K\"ahler metric which is
conformal to an Einstein metric is a constant nonzero multiple of the square root of the
conformal factor (cf. \cite[Proposition 6.5]{dr-ma1}). That square root is $\t$ for $g$,
and $\t/(\t-c)$ for $\wht{g}=g/(\t-c)^2$. Thus, writing $s=a\,\t$, $\wht{s}=b\,\t/(\t-c)$ for nonzero constants
$a$, $b$ and $\om$ for the K\"ahler form of $g$ we get
\be\lb{eta-end}
\eta_\gamma(\del M)=\fr 1{288\pi^2}\int_M \left[\left(b^2(\t-c)^{-6}-a^2\right)\t^2\right]\om^2/2-\sig(M).
\end{equation}

We now apply the Duistermaat-Heckman Theorem \cite{dh} to the integral in \Ref{eta-end}. Namely, the Killing potential $\t$ is a
moment map for the symplectic form $\omega$ associated with $g$, with respect to the circle action guaranteed by Theorem \ref{bihol-iso}.
The value $\t=0$ is regular for this moment map, because $|\nab\t|_{\t=0}\ne 0$, by \Ref{bdps}. The theorem states that the
Duistermaat-Heckman measure, i.e. the moment map push-forward of the measure associated with the symplectic form on $M$,
is the function $[\om_0+t\,\Omega][N]$, with $\om_0$ denoting the induced symplectic form on the associated reduced space
$N=\t^{-1}(0)/S^1$ and $\Omega$ denoting the curvature of $L$ (See \cite[Excercise $2.18$]{g}). We have
\begin{multline}
\eta_\gamma(\del M)=\fr 1{288\pi^2}\int_{0}^{\t_0} \left[\left(b^2(t-c)^{-6}-a^2\right)t^2\right]
\left(\left[\om_0+t\,\Omega\right]\left[N\right]\right)\,dt-\sig(M)\\
=\fr 1{288\pi^2}\int_{0}^{\t_0} \left[\left(b^2(t-c)^{-6}-a^2\right)t^2\right]
\left(
\left[\om_0+t\,p\,\om^h\right]\left[N\right]\right)\,dt-\sig(M)\\
=\fr 1{288\pi^2}Vol_h(N)\int_{0}^{\t_0} \left[\left(b^2(t-c)^{-6}-a^2\right)t^2\right]\left[l+p\,t\,\right]\,dt-\sig(M).
\end{multline}
The second line follows from \Ref{curv-L}.
The third follows from cases ii] and iii] of \Ref{met}, as $\om_0$ is the solution of the equation $\pi^*\om_0=\iota^*\om$,
valid on $\del M$ with $\iota$ being the inclusion $\del M \hookrightarrow M$ and $\pi$ the projection of $\del M$ to the
reduced space $\del M/S^1$.
\end{proof}

\begin{Remark}
The value of the signature $\sigma(M)$ in the above formula is in fact $0$ or $-1$ in the case where $M$ is homeomorphic to
a disk bundle $S$ over a Riemann surface $N$. Contractibility of
the disk implies that $H_2(S)=H_2(N)=\mathbf{Z}$, so that the intersection form is a number, namely the Euler
number of the bundle $S$. The group $H_2(S)$ is generated by a section, and the intersection matrix is determined
by its nonpositive self-intersection number.
\end{Remark}

\subsection{The ambitoric case}

Ambihermitian structures were defined in \cite{acg}. These consist, on a $4$-manifold $M$, of
a conformal class $c$, with two $c$-orthogonal complex structures $J$, $\wht{J}$ inducing opposite
orientations. A metric in such a conformal class is called ambihermitian. Suppose such a metric
$g_E$ is Einstein, and degeneracy of the Weyl tensor components $W_+$ and $W_-$ holds, in the sense that as operators on
selfdual  (respectfully antiselfdual) $2$-forms, at least two of their three eigenvalues coincide.
Then in a neighborhood of each point the structure is ambi-K\"ahler, i.e. there exist two metrics
$g$, $\wht{g}$ in the conformal class $c$ such that $(g,J)$, $(\wht{g},\wht{J})$ are K\"ahler \cite[Proposition 1]{acg}.
The three metrics satisfy $g_E = s^{-2} g = \wht{s}^{\,-2}\wht{g}$ on an open and dense set (away from the zero sets of
$s$ or $\wht{s}$), where $s$ ($\wht{s}$) is the scalar curvature of $g$ ($\wht{g}$). Furthermore,
$J\nab s$  and $\wht{J}\wht{\nab} \wht{s}$ are
Killing vector fields with respect to both $g$ and $\wht{g}$ \cite[Proposition 7]{acg}, so that
$\wht{J}\wht{\nab} \wht{s}=J\nab f$ and $J\nab s=\wht{J}\wht{\nab}\wht{f}$ for two functions $f$, $\wht{f}$.
Generically, these vector fields span a $2$-dimensional space of Poisson-commuting Killing vector fields,
and the whole structure is called ambitoric. We will also assume that $(s,f)$ and $(\wht{f},\wht{s})$ are
each moment maps for a toric group action (in the case of a manifold without boundary, this follows
from the assumption of having an ambitoric structure). We now show that in this case, a formula also
exists for the eta invariant of the boundary.
\begin{Theorem}\lb{ambi}
Let $M^4$ be a compact four-manifold having a nonempty boundary with a complete Einstein metric $g_E$
whose conformal class admits an ambitoric structure consisting of the metrics $g$, $\wht{g}$, with all metrics
defined on the interior. If the vector fields $J\nab s$  and $\wht{J}\wht{\nab} \wht{s}$ span an ambitoric
structure and each gives rise to a toric structure as described above,
then the eta invariant of the boundary metric $\gamma=g|_{\del M}$ is given by
$$\eta_\gamma(\del M)=\fr 1{288\pi^2}\int_{\wht{P}} \wht{y}^{\,2}\,v_{euc}-\fr 1{288\pi^2}\int_P x^2\,v_{euc}-\sig(M),$$
where $P$, $\wht{P}$ are the images of the moment maps $(s, f)$, $(\wht{f}, \wht{s})$ respectively.
Here $x$,
$\wht{y}$ are the functions on these images of the moment maps whose pullbacks are
$s$ and
$\wht{s}$, respectively, and $v_{euc}$ represents integration with respect to Lebesgue measure.
\end{Theorem}
\begin{proof}
Suppose $f$ is the function such that $g =  f^2\,\wht{g}$. We begin by noticing that
the relation generalizing \Ref{hat} is $|W_-|^2=f^{-4}|\wht{W}_+|_{\wht{g}}^2$, so that by \Ref{hitchin}
$$\eta_\gamma(\del M)=\fr 1{12\pi^2}\int_M \left[f^{-4}|\wht{W}_+|_{\wht{g}}^2-|W_+|^2\right]\text{vol}_g-\sig(M).$$
Next,   $\text{vol}_g=f^4\text{vol}_{\wht{g}}$, so that in the first line below, $f$ cancels to give:
\lb{ambi-comp}
\begin{eqnarray*}
\eta_\gamma(\del M)
&=&\fr 1{12\pi^2}\int_M |\wht{W}_+|_{\wht{g}}^2\,\text{vol}_{\wht{g}}-\fr 1{12\pi^2}\int_M|W_+|^2\text{vol}_g-\sig(M)\\
&=&\fr 1{288\pi^2}\int_M \wht{s}^{\ 2}\text{vol}_{\wht{g}}-\fr 1{288\pi^2}\int_M s^2\text{vol}_g-\sig(M)\\
&=&\fr 1{288\pi^2}\int_{\wht{P}} \wht{y}^{\,2}\,v_{euc}-\fr 1{288\pi^2}\int_P x^2\,v_{euc}-\sig(M).
\end{eqnarray*}
The last line follows from the fact that the Duistermaat-Heckman measure in the toric case equals
Lebesgue measure on the image of the moment map.
\end{proof}


\vspace{.3in}
\begin{center}
\textbf{Acknowledgements}
\end{center}
The author thanks the anonymous referee for a comment that led to a corrected statement
and improved proof outline of Theorem \ref{bihol-iso}.


\end{document}